\documentclass{amsart}
\usepackage{amsmath}
\usepackage{amssymb}

\newtheorem{thm}{Theorem}[section]
\newtheorem{cor}[thm]{Corollary}
\newtheorem{lem}[thm]{Lemma}

\theoremstyle{definition}
\newtheorem{defn}[thm]{Definition}
\newtheorem{rem}[thm]{Remark}
\newtheorem{exam}[thm]{Example}
\numberwithin{equation}{section}

\begin{document}
\title[Coupled fixed point results]{Some coupled fixed point results on partial metric spaces}
\author[H. Aydi
]
{
Hassen Aydi
}
\maketitle

\begin{abstract}
In this paper we give some coupled fixed point results for mappings satisfying different contractive
conditions on complete partial metric spaces.
\end{abstract}

\vspace{.25cm}\noindent\textbf{2000 Mathematics Subject
Classification.} 47H10, 54H25.

\vspace{.25cm}\noindent \textbf{Key Words and Phrases}: Coupled fixed
point, partial metric space, contractive condition

\section{\textbf{INTRODUCTION AND PRELIMINARIES}}
For a given partially ordered set $X$, Bhaskar and Lakshmikantham in \cite{BL} introduced the
concept of coupled fixed point of a mapping $F : X\times X \rightarrow X$ . Later in \cite{CL}, Lakshmikantham
and C\'ir\'ic investigated some more coupled fixed point theorems in partially ordered sets. The
following is the corresponding definition of a coupled fixed point.
\begin{defn}\cite{BL}
\label{def1}
An element $(x, y)\in X\times X$ is said to be a
coupled fixed point of the mapping $F : X\times X \rightarrow X$ if $F(x, y)=x$ and $F(y, x)=y$.
\end{defn}
F. Sabetghadam et al. \cite{SMS} obtained the following
\begin{thm}
Let $(X, d)$ be a complete cone metric space. Suppose that the mapping $F : X\times X \rightarrow X$ satisfies the following contractive condition for all $x$, $y, u, v\in X$
\[
d(F(x, y), F(u, v))\leq k d(x, u)+ld(y, v)
\]
where $k$, $l$ are nonnegative constants with $k+l < 1$. Then $F$ has a unique coupled fixed point.
\end{thm}
 In this paper, we give the analogous of
 this result (and some others in \cite{SMS}) on partial metric spaces, and we
establish some coupled fixed point results.\\
The concept of partial metric space $(X, p)$ was introduced by Matthews in 1994. In
such spaces, the distance of a point in the self may not be zero.
First, we start with some preliminaries definitions on the partial
metric spaces \cite{ASS, A, Ma2, O1, O2, OV, Ro1, Ro2, S1, V}
\begin{defn}\rm{(\cite{Ma2,O1,O2})}
\label{partial} A partial metric on a nonempty set $X$ is a
function $p : X \times X\longrightarrow \mathbb{R}^+$  such that for all
$x, y, z \in
X$:\\
\rm{(p1)} $x =y\Longleftrightarrow p(x, x) = p(x, y) = p(y, y)$,\\
\rm{(p2)} $p(x, x) \leq  p(x,y)$,\\
\rm{(p3)} $p(x, y)= p(y,x),$\\
\rm{(p4)} $p(x,y)\leq p(x,z)+ p(z,y) - p(z,z)$.\\
 A partial metric
space is a pair $(X, p)$ such that $X$ is a nonempty set and $p$
is a partial metric on $X$.
\begin{rem}
\label{rem0}
It is clear that, if $p(x,y)=0$, then from (p1), (p2) and (p3), $x=y$. But if $x
=y$, $p(x,y)$ may not be $0$.
\end{rem}
\end{defn}
If $p$ is a partial metric on $X$, then the function $p^s : X
\times X\longrightarrow R_+$ given by
\[
p^s(x, y) = 2p(x, y)- p(x, x)-p(y, y),
\]
 is a metric on $X$.
\begin{defn}\rm{(\cite{Ma2,O1,O2})}
\label{c-c}
Let $(X, p)$ be a partial metric space. Then: \\
\rm{(i)} a sequence $\{x_n\}$ in a partial metric space $(X, p)$
converges to a point $x\in X$ if and only if $p(x, x) =
\displaystyle\lim_{n\longrightarrow +\infty} p(x,
x_n)$;\\
\rm{(ii)} A sequence $\{x_n\}$ in a partial metric space $(X, p)$
is called a Cauchy sequence if there
exists (and is finite) $\displaystyle\lim_{n,m\longrightarrow +\infty} p(x_n, x_m)$. \\
\rm{(iii)} A partial metric space $(X, p)$ is said to be complete
if every Cauchy sequence $\{x_n\}$ in $X$ converges to a point
$x\in X$, that is $p(x, x)
=\displaystyle\lim_{n,m\longrightarrow +\infty} p(x_n, x_m)$.\\
\end{defn}
\begin{lem}\rm{(\cite{Ma2,O1,OV})}
\label{lem1}
 Let $(X, p)$ be a partial metric space.\\
(a) $\{x_n\}$ is a Cauchy sequence in $(X, p)$ if and only if it
is a Cauchy sequence in the metric space $(X, p^s)$. \\
(b) A partial metric space $(X, p)$ is complete if and only if the
metric space $(X, p^s)$ is complete. Furthermore,
$\displaystyle\lim_{n\longrightarrow +\infty} p^s(x_n, x) = 0$ if
and only if
\[
p(x, x) = \displaystyle\lim_{n\longrightarrow +\infty} p(x_n, x)
=\displaystyle\lim_{n,m\longrightarrow +\infty} p(x_n, x_m).
\]
\end{lem}
\section{\textbf{MAIN RESULTS}}
Our first main result is the following
 \begin{thm}
 \label{main-thm1}
 Let $(X, p)$ be a complete partial metric space. Suppose that the mapping $F : X\times X \rightarrow X$ satisfies the following contractive condition for all $x$, $y, u, v\in X$
\begin{equation}
\label{contraction}
p(F(x, y), F(u, v))\leq k p(x, u)+lp(y, v)
\end{equation}
where $k$, $l$ are nonnegative constants with $k+l < 1$. Then $F$ has a unique coupled fixed point.
\end{thm}
\noindent \textbf{Proof.} Choose $x_0$, $y_0\in X$ and set $x_1= F(x_0, y_0)$ and $y_1=F(y_0, x_0)$. Repeating this process, set $x_{n+1}= F(x_n, y_n)$ and $y_{n+1}=F(y_n, x_n)$. Then by (\ref{contraction}), we have
\begin{equation}
\label{eq1}
\begin{split}
p(x_n,x_{n+1})=&p(F(x_{n-1}, y_{n-1}),F(x_n, y_n))\\
\leq & kp(x_{n-1},x_{n})+l p(y_{n-1},y_{n}),
\end{split}
\end{equation}
and similarly
\begin{equation}
\label{eq2}
\begin{split}
p(y_n,y_{n+1})=&p(F(y_{n-1}, x_{n-1}),F(y_n, x_n))\\
\leq & kp(y_{n-1},y_{n})+l p(x_{n-1},x_{n}).
\end{split}
\end{equation}
Therefore, by letting
\begin{equation}
\label{eq3}
d_n=p(x_n,x_{n+1})+p(y_n,y_{n+1}),
\end{equation}
we have
\begin{equation}
\label{eq5}
\begin{split}
d_n=&p(x_n,x_{n+1})+p(y_n,y_{n+1})\\
\leq &kp(x_{n-1},x_{n})+l p(y_{n-1},y_{n})+kp(y_{n-1},y_{n})+l p(x_{n-1},x_{n})\\
=& (k+l)[p(y_{n-1},y_{n})+p(x_{n-1},x_{n})]\\
=& (k+l)d_{n-1}.
\end{split}
\end{equation}
Consequently, if we set $\delta= k+l$ then for each $n\in \mathbb{N}$ we have
\begin{equation}
\label{eq6}
d_n\leq \delta d_{n-1}\leq \delta^2 d_{n-2}\leq...\leq \delta^n d_0.
\end{equation}
If $d_0=0$ then $p(x_0,x_{1})+p(y_0,y_{1})=0$. Hence, from Remark \ref{rem0}, we get $x_0=x_1=F(x_0,y_0)$ and $y_0=y_1=F(y_0,x_0)$, meaning that $(x_0, y_0)$ is a coupled fixed point of $F$. Now, let $d_0 > 0$. For each $n\geq m$ we have in view of the condition $(p4)$
\[
\begin{split}
p(x_n,x_m)\leq& p(x_{n},x_{n-1})+p(x_{n-1},x_{n-2})-p(x_{n-1},x_{n-1})\\
+ & p(x_{n-2},x_{n-3})+p(x_{n-3},x_{n-4})-p(x_{n-3},x_{n-3})+\\
+&...+p(x_{m+2},x_{m+1})+p(x_{m+1},x_{m})-p(x_{m+1},x_{m+1})\\
\leq& p(x_{n},x_{n-1})+p(x_{n-1},x_{n-2})+...+p(x_{m+1},x_{m}).
\end{split}
\]
Similarly, we have
\[
p(y_n,y_m)\leq p(y_{n},y_{n-1})+p(y_{n-1},y_{n-2})+...+p(y_{m+1},y_{m}).
\]
Thus,
\begin{equation}
\label{eq7}
\begin{split}
p(x_n,x_m)+p(y_n,y_m)\leq &d_{n-1}+d_{n-2}+...+d_m\\
\leq & (\delta^{n-1}+\delta^{n-2}+...+\delta^{m})d_0\\
\leq & \frac{\delta^{m}}{1-\delta} d_0.
\end{split}
\end{equation}
By definition of $p^s$, we have $p^s(x,y)\leq 2 p(x,y)$, so for any $n\geq m$
\begin{equation}
\label{eq8}
p^s(x_n,x_m)+p^s(y_n,y_m)\leq 2p(x_n,x_m)+2p(y_n,y_m)\leq 2\frac{\delta^{m}}{1-\delta} d_0.
\end{equation}
which implies that $\{x_n\}$ and $\{y_n\}$ are Cauchy sequences in $(X,p^s)$ because of $0\leq\delta=k+l<1$. Since the partial metric space $(X,p)$ is complete, hence thanks to Lemma \ref{lem1}, the metric space $(X,p^s)$ is complete, so there exist $u^*$, $v^*\in X$ such
that
\begin{equation}
\label{eq9}
\lim_{n\rightarrow +\infty} p^s(x_n,u^*)=\lim_{n\rightarrow +\infty}p^s(y_n,v^*)=0.
\end{equation}
Again, from Lemma \ref{lem1}, we get
\[
p(u^*,u^*)=\lim_{n\rightarrow +\infty} p(x_n,u^*)=\lim_{n\rightarrow +\infty}p(x_n,x_n),
\]
and
\[
p(v^*,v^*)=\lim_{n\rightarrow +\infty} p(y_n,v^*)=\lim_{n\rightarrow +\infty}p(y_n,y_n).
\]
But, from condition $(p2)$ and (\ref{eq6}),
\[
p(x_n,x_n)\leq p(x_n,x_{n+1})\leq d_n\leq \delta^n d_0,
\]
so since $\delta\in [0,1[$, hence letting $n\rightarrow +\infty$, we get $\displaystyle\lim_{n\rightarrow +\infty} p(x_n,x_n)=0$. It follows that
\begin{equation}
\label{eq10}
p(u^*,u^*)=\lim_{n\rightarrow +\infty} p(x_n,u^*)=\lim_{n\rightarrow +\infty}p(x_n,x_n)=0.
\end{equation}
Similarly, we get
\begin{equation}
\label{eq11}
p(v^*,v^*)=\lim_{n\rightarrow +\infty} p(y_n,v^*)=\lim_{n\rightarrow +\infty}p(y_n,y_n)=0.
\end{equation}
Therefore, we have using (\ref{contraction})
\[
\begin{split}
p(F(u^*,v^*),u^*)\leq& p(F(u^*,v^*),x_{n+1})+p(x_{n+1},u^*)-p(x_{n+1},x_{n+1}),\quad\mbox{By (p4)}\\
\leq& p(F(u^*,v^*),F(x_{n},y_n))+p(x_{n+1},u^*)\\
\leq &k p(x_n,u^*)+l p(y_n,v^*)+p(x_{n+1},u^*),
\end{split}
\]
and letting $n\rightarrow +\infty$, then from (\ref{eq10}) and (\ref{eq11}), we obtain $p(F(u^*,v^*),u^*))=0$, so $F(u^*,v^*)=u^*$.
Similarly, we have $F(v^*,u^*)=v^*$, meaning that $(u^*,v^*)$ is a coupled fixed point of $F$.\\
Now, if $(u',v')$ is another coupled fixed point of $F$, then
\[
p(u',u^*)=p(F(u',v'),F(u^*,v^*))\leq k p(u',u^*)+l p(v',v^*)
\]
\[
p(v',v^*)=p(F(v',u'),F(v^*,u^*))\leq k p(v',v^*)+l p(u',u^*).
\]
It follows that
\[
p(u',u^*)+p(v',v^*)\leq (k+l)[p(u',u^*)+p(v',v^*)].
\]
In view of $k+l<1$, this implies that $p(u',u^*)+p(v',v^*)=0$,  so $u^*=u'$ and $v^*=v'$. The proof of Theorem \ref{main-thm1} is completed.\\

It is worth noting that when the constants in Theorem  \ref{main-thm1} are equal we have the
following Corollary
\begin{cor}
\label{cor1}
Let $(X, p)$ be a complete partial metric space. Suppose that the mapping $F : X\times X \rightarrow X$ satisfies the following contractive condition for all $x$, $y, u, v\in X$
\begin{equation}
\label{contraction1}
p(F(x, y), F(u, v))\leq \frac{k}{2}( p(x, u)+p(y, v))
\end{equation}
where  $0\leq k< 1$. Then, $F$ has a unique coupled fixed point.
\end{cor}
\begin{exam}
\label{ex1}
Let $X=[0,+\infty[$ endowed with the usual partial metric $p$
defined by $p:X\times X\rightarrow [0,+\infty[$ with
$p(x,y)=\max\{x,y\}$. The partial metric space
 $(X,p)$ is complete because $(X,p^s)$ is complete. Indeed,
 for any $x,y\in X$,
 \[
 \begin{split}
 p^s(x,y)=2p(x,y)-p(x,x)-p(y,y)=&2\max\{x,y\}-(x+y)\\
 =& |x-y|,
 \end{split}
 \]
Thus, $(X, p^s)$ is the Euclidean
metric space which is complete. Consider the mapping
 $F:X\times X\rightarrow X$ defined by $F(x,y)=\frac{x+y}{6}$. For any $x,y,u,v\in X$, we have
 \[
 p(F(x,y),F(u,v))=\frac{1}{6}\max\{x+y, u+v\}\leq \frac{1}{6}[\max\{x,u\}+\max\{y,v\}]=\frac{1}{6} [p(x,u)+p(y,v)],
 \]
which is the contractive condition (\ref{contraction1}) for $k=\frac{1}{3}$. Therefore, by Corollary \ref{cor1}, $F$ has a unique coupled fixed point, which is $(0, 0)$.
Note that if the mapping $F : X\times X \rightarrow X$ is given by $F(x,y)=\frac{x+y}{2}$, then $F$ satisfies the
contractive condition (\ref{contraction1}) for $k=1$, that is,
\[
 p(F(x,y),F(u,v))=\frac{1}{2}\max\{x+y, u+v\}\leq \frac{1}{2}[\max\{x,u\}+\max\{y,v\}]=\frac{1}{2} [p(x,u)+p(y,v)],
 \]
In this case, $(0, 0)$ and $(1, 1)$ are both coupled fixed points of $F$ and hence the coupled fixed
point of $F$ is not unique. This shows that the condition $k < 1$ in Corollary \ref{cor1}, and hence
$k+l < 1$ in Theorem \ref{main-thm1} can not be omitted in the
statement of the aforesaid results.
\end{exam}
\begin{thm}
\label{main-thm2}
 Let $(X, p)$ be a complete partial metric space. Suppose that the mapping $F : X\times X \rightarrow X$ satisfies the following contractive condition for all $x$, $y, u, v\in X$
\begin{equation}
\label{contraction2}
p(F(x, y), F(u, v))\leq k p(F(x,y),x)+lp(F(u,v),u)
\end{equation}
where $k$, $l$ are nonnegative constants with $k+l < 1$. Then $F$ has a unique coupled fixed point.
\end{thm}
\noindent \textbf{Proof.}
We take the same sequences $\{x_n\}$ and $\{y_n\}$ given in the proof of Theorem \ref{main-thm1}  by
\[
x_{n+1}=F(x_n,y_n),\quad y_{n+1}=F(y_n,x_n)\quad\mbox{for any}\quad n\in\mathbb{N}.
\]
Applying (\ref{contraction2}), we get
\begin{equation}
\label{eq12}
p(x_n,x_{n+1})\leq \delta p(x_{n-1},x_{n})
\end{equation}
\begin{equation}
\label{eq13}
p(y_n,y_{n+1})\leq \delta p(y_{n-1},y_{n}),
\end{equation}
where $\delta =\frac{k}{1-l}$. By definition of $p^s$, we have
\begin{equation}
\label{eq14}
p^s(x_n,x_{n+1})\leq 2p(x_n,x_{n+1}) \leq 2\delta^n p(x_{1},x_{0})
\end{equation}
\begin{equation}
\label{eq15}
p^s(y_n,y_{n+1})\leq 2 p(y_n,y_{n+1})\leq 2\delta^n p(y_{1},y_{0}).
\end{equation}
 Since $k+l<1$, hence $\delta<1$, so the sequences $\{x_n\}$ and $\{y_n\}$ are Cauchy sequences in the metric space $(X,p^s)$. The partial metric space $(X,p)$ is complete, hence from Lemma \ref{lem1}, $(X,p^s)$ is complete, so there
exist
 $u^*$, $v^*\in X$ such
that
\begin{equation}
\label{eq16}
\lim_{n\rightarrow +\infty} p^s(x_n,u^*)=\lim_{n\rightarrow +\infty}p^s(y_n,v^*)=0.
\end{equation}
From Lemma \ref{lem1}, we get
\[
p(u^*,u^*)=\lim_{n\rightarrow +\infty} p(x_n,u^*)=\lim_{n\rightarrow +\infty}p(x_n,x_n),
\]
and
\[
p(v^*,v^*)=\lim_{n\rightarrow +\infty} p(y_n,v^*)=\lim_{n\rightarrow +\infty}p(y_n,y_n).
\]
By the condition (p2) and (\ref{eq12}), we have
\[
p(x_n,x_{n})\leq p(x_n,x_{n+1}) \leq \delta^n p(x_{1},x_{0}),
\]
so $\displaystyle\lim_{n\rightarrow +\infty}p(x_n,x_n)=0$. It follows that
\begin{equation}
\label{eq17}
p(u^*,u^*)=\lim_{n\rightarrow +\infty} p(x_n,u^*)=\lim_{n\rightarrow +\infty}p(x_n,x_n)=0.
\end{equation}
Similarly, we find
\begin{equation}
\label{eq18}
p(v^*,v^*)=\lim_{n\rightarrow +\infty} p(y_n,v^*)=\lim_{n\rightarrow +\infty}p(y_n,y_n)=0.
\end{equation}
Therefore, by (\ref{contraction2})
\[
\begin{split}
p(F(u^*,v^*),u^*)\leq& p(F(u^*,v^*),x_{n+1})+p(x_{n+1},u^*)\\
=& p(F(u^*,v^*),F(x_{n},y_n))+p(x_{n+1},u^*)\\
\leq &k p(F(u^*,v^*),u^*)+l p(F(x_n,y_n),x_n)+p(x_{n+1},u^*)\\
=&k p(F(u^*,v^*),u^*)+l p(x_{n+1},x_n)+p(x_{n+1},u^*)
\end{split}
\]
and letting $n\rightarrow +\infty$, then from (\ref{eq14})-(\ref{eq17}), we obtain
\[
p(F(u^*,v^*),u^*)\leq k p(F(u^*,v^*),u^*).
\]
From the preceding inequality we can
deduce a contradiction if we assume that $p(F(u^*, v^*), u^*)\neq 0$, because
in that case we conclude that $1\leq k$ and now this inequality is, in fact,
a contradiction, so $p(F(u^*,v^*),u^*)=0$, that is, $F(u^*,v^*)=u^*$.
Similarly, we have $F(v^*,u^*)=v^*$, meaning that $(u^*,v^*)$ is a coupled fixed point of $F$.
Now, if $(u',v')$ is another coupled fixed point of $F$, then in view of (\ref{contraction2})
\[
\begin{split}
p(u',u^*)=& p(F(u',v'),F(u^*,v^*))\\
\leq& k p(F(u',v'),u')+l p(F(u^*,v^*),u^*)\\
=& k p(u',u')+l p(u^*,u^*)\\
\leq & k p(u',u^*)+l p(u',u^*)=(k+l) p(u',u^*),\quad\mbox{using (p2)}
\end{split}
\]
that is $p(u',u^*)=0$ since $(k+l)<1$.
It follows that $u^*=u'$. Similarly, we can have $v^*=v'$, and the proof of Theorem \ref{main-thm2} is completed.\\
\begin{thm}
\label{main-thm3}
 Let $(X, p)$ be a complete partial metric space. Suppose that the mapping $F : X\times X \rightarrow X$ satisfies the following contractive condition for all $x$, $y, u, v\in X$
\begin{equation}
\label{contraction3}
p(F(x, y), F(u, v))\leq k p(F(x,y),u)+lp(F(u,v),x)
\end{equation}
where $k$, $l$ are nonnegative constants with $k+2l<1$. Then $F$ has a unique coupled fixed point.
\end{thm}
\noindent \textbf{Proof.}
Since, $k+2l<1$, hence $k+l<1$, and as a consequence the proof of the uniqueness in this
Theorem is as trivial as in the other results. To prove the existence of the fixed point, choose the sequences
$\{x_n\}$ and $\{y_n\}$ like in the proof of Theorem \ref{main-thm1}, that is
\[
x_{n+1}=F(x_n,y_n),\quad y_{n+1}=F(y_n,x_n)\quad\mbox{for any}\quad n\in\mathbb{N}.
\]
Applying again (\ref{contraction3}), we have
\[
\begin{split}
p(x_n,x_{n+1})=&p(F(x_{n-1}, y_{n-1}),F(x_n, y_n))\\
\leq & kp(F(x_{n-1},y_{n-1}),x_{n})+l p(F(x_{n},y_{n}),x_{n-1})\\
=& kp(x_n,x_n)+l p(x_{n+1},x_{n-1})\\
\leq &k p(x_{n+1},x_{n})+l p(x_{n+1},x_{n-1})],\quad\mbox{by (p2)}\\
\leq& k p(x_{n+1},x_{n})+l p(x_{n+1},x_{n})+l p(x_{n},x_{n-1})-lp(x_n,x_n),\quad\mbox{using (p4)} \\
\leq& (k+l) p(x_{n},x_{n+1})+l p(x_{n-1},x_{n}).
\end{split}
\]
It follows that for any $n\in \mathbb{N}^*$
\[
p(x_n,x_{n+1})\leq \frac{l}{1-l-k}p(x_{n-1},x_{n}).
\]
Let us take $\delta=\frac{l}{1-l-k}$.  Hence, we deduce
\begin{equation}
\label{eq19}
p^s(x_n,x_{n+1})\leq 2 p(x_n,x_{n+1})\leq 2 \delta^n p(x_{0},x_{1}).
\end{equation}
Under the condition $0\leq k+2l<1$, we get $0\leq \delta<1$.
From this fact we immediately obtain that $\{x_n\}$ is Cauchy in the complete metric
space $(X, p^s)$. Of course, similar arguments apply to the
case of the sequence $\{y_n\}$ in order to prove that
\begin{equation}
\label{eq20}
p^s(y_n,y_{n+1})\leq 2 p(y_n,y_{n+1})\leq 2 \delta^n p(y_{0},y_{1}),
\end{equation}
and, thus, that the sequence $\{y_n\}$ is Cauchy in $(X, p^s)$.
Therefore, there exist
 $u^*$, $v^*\in X$ such
that
\begin{equation}
\label{eq21}
\lim_{n\rightarrow +\infty} p^s(x_n,u^*)=\lim_{n\rightarrow +\infty}p^s(y_n,v^*)=0.
\end{equation}
Thanks to Lemma \ref{lem1}, we have
\[
\lim_{n\rightarrow +\infty} p(x_n,u^*)=\lim_{n\rightarrow +\infty}p(x_n,x_n)=p(u^*,u^*),
\]
and
\[
\lim_{n\rightarrow +\infty} p(y_n,v^*)=\lim_{n\rightarrow +\infty}p(y_n,y_n)=p(v^*,v^*).
\]
The condition $(p2)$ together with (\ref{eq19}) yield that
\[
p(x_n,x_n)\leq p(x_n,x_{n+1})\leq \delta^n p(x_{0},x_{1}),
\]
hence letting $n\rightarrow +\infty$, we get $\displaystyle\lim_{n\rightarrow +\infty}p(x_n,x_n)=0$. It follows that
\begin{equation}
\label{eq22}
p(u^*,u^*)=\lim_{n\rightarrow +\infty} p(x_n,u^*)=\lim_{n\rightarrow +\infty}p(x_n,x_n)=0.
\end{equation}
Similarly, we have
\begin{equation}
\label{eq23}
p(v^*,v^*)=\lim_{n\rightarrow +\infty} p(y_n,v^*)=\lim_{n\rightarrow +\infty}p(y_n,y_n)=0.
\end{equation}
Therefore, we have using (\ref{contraction3})
\[
\begin{split}
p(F(u^*,v^*),u^*)\leq& p(F(u^*,v^*),x_{n+1})+p(x_{n+1},u^*)\\
=& p(F(u^*,v^*),F(x_{n},y_n))+p(x_{n+1},u^*)\\
\leq &k p(F(u^*,v^*),x_n)+l p(F(x_n,y_n),u^*)+p(x_{n+1},u^*)\\
=&k p(F(u^*,v^*),x_n)+l p(x_{n+1},u^*)+p(x_{n+1},u^*)\\
\leq &k p(F(u^*,v^*),u^*)+k p(u^*,x_n)+l p(x_{n+1},u^*)+p(x_{n+1},u^*),\quad\mbox{using p(4)}.
\end{split}
\]
Letting $n\rightarrow +\infty$ yields, using (\ref{eq22})
\[
p(F(u^*,v^*),u^*)\leq k p(F(u^*,v^*),u^*),
\]
and since $k<1$, we have $p(F(u^*,v^*),u^*)=0$, that is $F(u^*,v^*)=u^*$. Similarly, thanks to (\ref{eq23}), we get $F(v^*,u^*)=v^*$, and hence  $(u^*,v^*)$ is a coupled fixed point of $F$.\\

When the constants in Theorems \ref{main-thm2} and \ref{main-thm3} are equal, we get the following corollaries
\begin{cor}
\label{cor2}
Let $(X, p)$ be a complete partial metric space. Suppose that the mapping $F : X\times X \rightarrow X$ satisfies the following contractive condition for all $x$, $y, u, v\in X$
\begin{equation}
\label{contraction4}
p(F(x, y), F(u, v))\leq \frac{k}{2}( p(F(x,y), x)+p(F(u,v), u))
\end{equation}
where  $0\leq k< 1$. Then, $F$ has a unique coupled fixed point.
\end{cor}
\begin{cor}
\label{cor3}
Let $(X, p)$ be a complete partial metric space. Suppose that the mapping $F : X\times X \rightarrow X$ satisfies the following contractive condition for all $x$, $y, u, v\in X$
\begin{equation}
\label{contraction5}
p(F(x, y), F(u, v))\leq \frac{k}{2}( p(F(x,y), u)+p(F(u,v), x))
\end{equation}
where  $0\leq k< \frac{2}{3}$. Then, $F$ has a unique coupled fixed point.
\end{cor}
\noindent \textbf{Proof.} The condition $0\leq k< \frac{2}{3}$ follows from the hypothesis on $k$ and $l$ given in Theorem \ref{main-thm3}.\\
\begin{rem}
\begin{itemize}
\item  Theorem \ref{main-thm1} extends the Theorem 2.2 of \cite{SMS} on the class of partial metric spaces.
\item  Theorem \ref{main-thm2} extends the Theorem 2.5 of \cite{SMS} on the class of partial metric spaces.
    \end{itemize}
    \end{rem}
\begin{rem}
Note that in Theorem \ref{main-thm2}, if the mapping $F : X\times  X\rightarrow X$ satisfies the contractive
condition (\ref{contraction2}) for all $x$, $y$, $u$, $v\in X$, then $F$ also satisfies the following contractive condition
\begin{equation}
\begin{split}
\label{contraction6}
p(F(x, y), F(u, v))=& p(F(u, v), F(x, y))\\
\leq& k p(F(u,v), u)+ l p(F(x,y), x)
\end{split}
\end{equation}
Consequently, by adding (\ref{contraction2}) and (\ref{contraction6}), $F$ also satisfies the following:
\begin{equation}
\label{contraction7}
p(F(x, y), F(u, v))\leq \frac{k+l}{2} p(F(u,v), u)+ \frac{k+l}{2} p(F(x,y), x)
\end{equation}
which is a contractive condition of the type (\ref{contraction4}) in Corollary \ref{cor2} with equal constants.
Therefore, one can also reduce the proof of general case (\ref{contraction2}) in Theorem \ref{main-thm2} to the special
case of equal constants. A similar argument is valid for the contractive conditions (\ref{contraction3}) in
Theorem \ref{main-thm3} and (\ref{contraction5}) in Corollary \ref{cor3}.

\end{rem}
\noindent{\bf Acknowledgment.} The author thanks the editor and the referees for their kind comments and suggestions to improve this paper.

\vspace{0.2cm}

\noindent Hassen Aydi:\newline Universit\'e de Monastir.\newline
Institut Sup\'erieur d'Informatique de Mahdia. Route de R\'ejiche,
Km 4, BP 35, Mahdia 5121, Tunisie.\newline Email-address:
hassen.aydi@isima.rnu.tn\newline

\end{document}